\documentclass[a4paper]{article}

\usepackage{anysize} 
\marginsize{3cm}{3cm}{1in}{1in} 
\usepackage{amsmath,amssymb,amsthm}
\usepackage{hyperref}
\hypersetup{
    pdfstartview={FitH},     									
    pdftitle={A Brief Survey of Higgs Bundles},             				
    pdfauthor={Ronald Alberto Z\'u\~niga-Rojas},     						
    pdfsubject={Vector bundles on curves and their moduli},           				
    pdfkeywords={Hodge Bundles, Moduli of Higgs Bundles, Stratifications, 
    Variations of Hodge Structures, Vector Bundles},                                            
    pdfnewwindow=true,        									
    colorlinks=true,          									
    linkcolor=blue,           									
    citecolor=red,            									
    urlcolor=blue             									
}
\usepackage[utf8]{inputenc}

\usepackage{enumerate}
\usepackage[T1]{fontenc}
\usepackage[dvips]{graphicx}

\usepackage{times}
\usepackage{multicol} 

\usepackage{makeidx}

\usepackage{xcolor}

\makeindex
\newcommand{\nada}[1]{}

\DeclareMathOperator{\End}{End}     
\DeclareMathOperator{\GCD}{GCD}     
\DeclareMathOperator{\rk}{rk}       

\newcommand{\Sg}{\Sigma}            

\newcommand{\bC}{\mathbb{C}}        
\newcommand{\bR}{\mathbb{R}}        
\newcommand{\sA}{\mathcal{A}}       
\newcommand{\sE}{\mathcal{E}}       
\newcommand{\sF}{\mathcal{F}}       
\newcommand{\gG}{\mathcal{G}}               
\newcommand{\sM}{\mathcal{M}}       
\newcommand{\sN}{\mathcal{N}}       

\renewcommand{\geq}{\geqslant}      
\renewcommand{\leq}{\leqslant}      
\newcommand{\ox}{\otimes}           
\renewcommand{\:}{\colon}           

\newcommand{\mathword}[1]{\quad\text{#1}\quad} 

\theoremstyle{plain}
\newtheorem{Th}{Theorem}[section]   
\newtheorem*{nonum-Th}{Theorem}     
\newtheorem{Prop}[Th]{Proposition}  
\newtheorem*{nonum-Cor}{Corollary}  

\theoremstyle{definition}
\newtheorem{Def}[Th]{Definition}    

\theoremstyle{remark}
\newtheorem{Rmk}[Th]{Remark}        
\newtheorem*{nonum-Rmk}{Remark}  

\DeclareRobustCommand{\QEDA}{\ifmmode
  \else \leavevmode\unskip\penalty9999 \hbox{}\nobreak\hfill \fi
  \quad\hbox{\qedasymbol}}
\newcommand{\qedasymbol}{$\boxminus$} 

\begin{document}

\begin{center}
\Large
\textsc{A Brief Survey of Higgs Bundles\footnote{Presented as Communication in SIMMAC 2018.}}

\bigskip
\normalsize
  March 05, 2019\\
  
\bigskip
\emph{Ronald A. Z\'u\~niga-Rojas}\footnote{ \scriptsize
  Supported by Universidad de Costa Rica through Escuela de Matem\'atica, specifically through CIMM (Centro de Investigaciones Matem\'aticas y Metamatem\'aticas), 
  Project 820-B8-224. This work is partly based on the Ph.D. Project~\cite{z-r0} called
  ``Homotopy Groups of the Moduli Space of Higgs Bundles'', supported
  by FEDER through Programa Operacional Factores de Competitividade-COMPETE,
  and also supported by FCT (Funda\c{c}\~ao para a Ci\^encia e a Tecnologia)
  through the projects PTDC/MAT-GEO/0675/2012 and PEst-C/MAT/UI0144/2013 
  with grant reference SFRH/BD/51174/2010.}\\[6pt]
\small Centro de Investigaciones Matem\'aticas y Metamatem\'aticas CIMM\\
\small Escuela de Matem\'atica, Universidad de Costa Rica, San Jos\'e 11501, Costa Rica\\
\small e-mail: \texttt{ronald.zunigarojas@ucr.ac.cr}
\end{center}

\title{A Brief Survey of Higgs Bundles\footnote{ \scriptsize Presented as Communication in SIMMAC 2018.}\\ 
{\Large Un estudio conciso de fibrados de Higgs}
}

\author{\sc Z\'u\~niga-Rojas, Ronald Alberto\\ CIMM, Universidad de Costa Rica}



\begin{abstract}

Considering a compact Riemann surface of genus greater than two, a Higgs bundle is a pair composed of a holomorphic bundle over the Riemann surface, joint with an auxiliar vector field, so-called Higgs field. This theory started around thirty years ago, with Hitchin's work, when he reduced the self-duality equations from dimension four to dimension two, and so, studied those equations over Riemann surfaces. Hitchin baptized those fields as \emph{Higgs fields} because in the context of physics and gauge theory, they describe similar particles to those described by the Higgs bosson. Later, Simpson used the name \emph{Higgs bundle} for a holomorphic bundle together with a Higgs field. Today, Higgs bundles are the subject of research in 
several areas such as non-abelian Hodge theory, Langlands, mirror symmetry, integrable 
systems, quantum field theory (QFT), among others.
The main purposes here are to introduce these objects, and to present a brief but complete construction of the moduli space of Higgs bundles, and some of its stratifications.

\begin{flushleft}

\emph{MSC classes}:  Primary \texttt{14H60}; Secondaries \texttt{14D07}, \texttt{55Q52}.
\end{flushleft}

\bigskip

\centerline{{\bf Resumen}}\medskip

Considerando una superficie compacta de Riemann de g\'enero mayor que dos, un fibrado de Higgs es un par compuesto por un fibrado holomorfo sobre la superficie de Riemann, junto con un campo vectorial auxiliar, llamado campo de Higgs.Esta teor\'ia inici\'o hace unos treinta a\~nos, con el trabajo de Hitchin, cuando \'el reduce las ecuaciones de autodualidad de dimensi\'on cuatro a dimensi\'on dos, y as\'i, estudiar esas ecuaciones sobre superficies de Riemann. Hitchin bautiz\'o esos campos como \emph{campos de Higgs} pues en el contexto de la f\'isica y de la teor\'ia de gauge, describen part\'iculas similares a las descritas por el boz\'on de Higgs. M\'as tarde, Simpson us\'o el nombre \emph{fibrado de Higgs} para un fibrado holomorfo junto con un campo de Higgs. Hoy, los fibrados de Higgs son objeto de investigaci\'on en varias \'areas tales como la teor\'ia de Hodge no abeliana, Langlands, simetr\'ia de espejo, sistemas integrables, teor\'ia cu\'antica de campos (QFT), entre otros.Los prop\'ositos principales aqu\'i son introducir estos objetos y presentar una breve pero completa construcci\'on del espacio m\'oduli de los fibrados de Higgs y algunas de sus estratificaciones.

\begin{flushleft}
\small
\emph{Palabras clave}:
Fibrados de Higgs, Fibrados de Hodge, Espacios M\'oduli, Triples Estables, Fibrados Vectoriales.
\end{flushleft}

\end{abstract}

\section*{Introduction} 
\label{sec:0-intro} 

Nineteen years ago, the Clay Institute of Cambridge, Massachusets, presented seven millenium challenging problems. Yang-Mills and Mass Gap is one of those challenging unsolved problems. Newton's classical mechanics laws stand to planets and celestial bodies as Quantum laws stand to elementary particles. Sixty five years ago, C.~N.~Yang and R.~L.~Mills~\cite{yami} introduced a remarkable framework to describe particles using structures that occur in geometry. Quantum Yang-Mills theory is nowadays the foundation of elementary particle physics theory and its predictions have been tested at many laboratories; nevertheless, its mathematical foundations are still a mystery. Yang-Mills theory depends on a quantum mechanical property, {\em the mass gap}: the quantum particles have positive masses, even though the classical waves travel at the speed of light. Establishing the existence of the Yang-Mills theory and a mass gap will require the introduction of fundamental new ideas both in physics and in mathematics. Thirty two years ago, Hitchin~\cite{hit2} reduced Yang-Mills self-duality equations from $\bR^4$ to $\bR^2$ imposing invariance translation in two directions, and considering an auxiliar vector field, so-called {\em Higgs field}. Hitchin's reduction of Yang-Mills self-duality equations in $\bR^2$ have the important property of conformal invariance, which allows them to be defined on Riemann surfaces.

Let $\Sigma = \Sigma_g$ be a compact (closed and connected) Riemann surface of genus $g > 2$.
Let $K = K_{\Sigma} = T^* \Sigma$ be the canonical line bundle over $\Sigma$ (its cotangent bundle).

\begin{Rmk} 
 Algebraically, $\Sigma_g$ is also a complete irreducible nonsingular algebraic curve over $\bC$:
 \[
 \dim_{\bC}(\Sigma_g) = 1. 
 \]
\end{Rmk}

The paper is organized as follows: 
in section \ref{sec:1.} we recall some basic facts about differential geometry, specifically about smooth vector bundles and holomorphic vector bundles; in section \ref{sec:2.}, we present the gauge group action over the set of connections and the induced action over the holomorphic structures; in section \ref{sec:3.}, we define the moduli space of Higgs bundles, and we present the two principal constructions: Hitchin construction in \ref{ssec:3.2.} and Simpson construction in \ref{ssec:3.3.}; finally, in \ref{sec:4.}, we present the most recent result of the author, published in \cite{z-r2} in terms of homotopy.

\section{Preliminary definitions} 
\label{sec:1.}

We shall recall some basic definitions we will need, like bundles and some of their invariants 
such as rank, degree and slope. We also present the definition of stability.

\begin{Def} 
 For any smooth vector bundle $\sE \to \Sg$, we denote the rank of~$\sE$ by $\rk(\sE) = r$ and the degree of $\sE$ by $\deg(\sE) = d$. Then, for any smooth complex bundle $\sE \to \Sg$ the \emph{slope} is defined to be
  \begin{equation}
    \mu(\sE) := \frac{\deg(\sE)}{\rk(\sE)} = \frac{d}{r}.
    \label{slope} 
  \end{equation} 
  \end{Def}

\begin{Def} 
 A smooth vector bundle $\sE \to \Sg$ is called \emph{semistable} if 
 $\mu(\sF) \leq \mu(\sE)$ for any $\sF$ such that
 $0 \subsetneq \sF \subseteq \sE$. Similarly, a vector bundle $\sE \to \Sg$ is
 called \emph{stable} if $\mu(\sF) < \mu(\sE)$ for any nonzero proper
 subbundle $0 \subsetneq \sF \subsetneq \sE$. Finally, $\sE$ is called 
 \emph{polystable} if it is the direct sum of stable subbundles,
 all of the same slope.
\end{Def}

Let $\sE\to \Sg$ be a complex smooth vector bundle with a hermitian metric on it.

\begin{Def} 
  A \emph{connection} $d_A$ on $\sE$ is a differential operator
$$
d_A: \Omega^{0}(\Sg,\sE)\to \Omega^{1}(\Sg,\sE)
$$
such that
$$
d_A(fs)=df\otimes s + fd_A s
$$
for any function 
$f \in C^{\infty}(\Sg)$ and any section 
$s \in \Omega^{0}(\Sg,\sE)$ where $\Omega^{n}(\Sg,\sE)$ 
is the set of smooth $n$-forms of 
$\Sg$ with values in $\sE$. Locally:
$$
d_A=d+A=d+Cdz+Bd\bar{z}
$$
where $A$ is a matrix of $1$-forms: $A_{ij}\in \Omega^{1}(\Sg,\sE)$, 
and $B, C$ are matrix valued functions depending 
on the hermitian metric on $\sE$.
\end{Def}

\begin{Def} 
When a connection $d_A$ is compatible with the hermitian metric on $\sE$, {\em i.e.}
\[
d \langle s,t\rangle = \langle d_A s,t\rangle + \langle s,d_A t\rangle 
\]
for the hermitian inner product $\langle \cdot, \cdot\rangle$ and for $s,t$ any couple 
of sections of $\sE$, $d_A$ is a \emph{unitary} connection. Denote $\sA(\sE)$ as the space of 
unitary connections on $\sE$, for a smooth bundle $\sE \to \Sg$. 
\end{Def}

\begin{Rmk} 
 Some authors call the matrix $A$ as a connection and call 
 $d_A=d+A$ as its correspoding covariant derivative. We abuse 
 notation and will not distinguish between them.
\end{Rmk}

\begin{Def} 
 The fundamental invariant of a connection is its \emph{curvature}:
$$
F_A:=d_{A}^{2}=d_A\circ d_A:\Omega^{0}(\Sg,\sE)\to \Omega^{2}(\Sg,\sE)
$$
where we are extending $d_A$ to $n$-forms in $\Omega^{n}(\Sg,\sE)$ in the obvious way. Locally:
\[
 F_A = dA + A^2.
\]
$F_A$ is $C^{\infty}(\Sg,\sE)$-linear and can be considered as a $2$-form on $\Sg$ with values in 
$\mathrm{End}(\sE):\ F_A\in \Omega^{2}(\Sg,\mathrm{End}(\sE))$, or locally as a matrix-valued $2$-form.
\end{Def}

\begin{Def} 
 If the curvature vanishes $F_A=0$, we say that the connection $d_A$ is \emph{flat}. 
 A flat connection gives a family of constant transition functions for $\sE$, which in 
 turn defines a representation of the fundamental group $\pi_1(\Sg)$ into $GL_r(\mathbb{C})$,
 and the image is in $U(r)$ if $d_A$ is unitary. Besides, from Chern-Weil theory, if $F_A = 0$, 
 then $\mathrm{deg}(\sE) = 0$.
\end{Def}

\section{Gauge Group Action} 
\label{sec:2.}

The gauge theory arises first in the context of physics, specifically in general relativity and 
classical~electromagnetism. Here, we present the mathematical formalism.

 \begin{Def} 
 A \emph{gauge transformation} is an automorphism of $\sE$. 
 Locally, a gauge transformation $g\in Aut(\sE)$ is a 
 $C^{\infty}$-function with values in $GL_r(\bC)$. A gauge 
 transformation $g$ is \emph{unitary} if $g$ preserves the 
 hermitian inner product. 

We denote by $\gG$ the group of unitary gauge transformations. This gauge group $\gG$ acts on 
$\sA(\sE)$ by conjugation:
$$
g \cdot d_A = g^{-1} d_A g\ 
\forall g \in \gG \mathword{and for} 
d_A \in \sA(\sE).
$$
 \end{Def}

\begin{Rmk} 
 Note that conjugation by a unitary gauge transformation takes 
 a unitary connection to a unitary connection.
 We denote by $\bar{\gG}$ the quotient of $\gG$ by the central subgroup $U(1)$,
 and by $\gG^{\bC}$ the complex gauge group. Besides, we denote by $B\gG$ and by
 $B\bar{\gG}$ the classifying spaces of $\gG$ and $\bar{\gG}$ respectively.
  There is a homotopical equivalence:
 $\gG \simeq \gG^{\bC}$.
\end{Rmk}

\begin{Def} 
 A \emph{holomorphic structure} on $\sE$ is a differential operator:
\[
  \bar{\partial}_{A}: \Omega^{0}(\Sg,\sE)\to \Omega^{0,1}(\Sg,\sE)
\]
such that $\bar{\partial}_{A}$ satisfies Liebniz rule:
\[
\bar{\partial}_{A} (fs)=\bar{\partial}f\otimes s + f\bar{\partial}_{A} s
\]
and \emph{the integrability condition}, which is precisely vanishing of 
the curvature 
\[
 F_A = 0.
\]
Here, $\bar{\partial}f=\frac{\partial f}{\partial \bar{z}} d\bar{z}$, 
and $\Omega^{p,q}(\Sg,\sE)$ is the space of smooth $(p,q)$-forms 
with values in $\sE$. Locally:
$$
\bar{\partial}_{A}=\bar{\partial}+A^{0,1} d\bar{z}
$$
where $A^{0,1}$ is a matrix valued function. The reader may consult the details in \cite{hit1}, \cite{kob} or \cite{lute}.
\end{Def}

\begin{Rmk} 
 Denote the space of holomorphic structures on $\sE$ by $\sA^{0,1}(\sE)$, and consider the collection
 $\sA^{0,1}(r,d) = \{\sA^{0,1}(\sE)\: \rk(\sE) = r \mathword{and} \deg(\sE) = d\}$. Here, the 
 subcollections  
 $$
 \sA_{ss}^{0,1}(r,d) = \{\sA^{0,1}(\sE)\: \sE \mathword{is semistable}\} 
 \subseteq \sA^{0,1}(r,d),
 $$
 
 $$
 \sA_{s}^{0,1}(r,d) = \{\sA^{0,1}(\sE)\: \sE \mathword{is stable}\} 
 \subseteq \sA^{0,1}(r,d),
 $$
  
 $$
 \sA_{ps}^{0,1}(r,d) = \{\sA^{0,1}(\sE)\: \sE \mathword{is polystable}\} \subseteq \sA^{0,1}(r,d),
 $$
 
 will be of particular interest for us.
\end{Rmk}

\section[M.S.~of~Higgs~bundles]{Moduli Space of Higgs Bundles} 
\label{sec:3.}

A central feature of Higgs bundles, is that they come in collections parametrized by the 
points of a quasi-projective variety: the Moduli Space of Higgs Bundles.

\subsection[M.S.~of~vector~bundles]{Moduli space of vector bundles} 
\label{ssec:3.1.}

The starting point to classify these collections, is the Moduli Space of Vector bundles. 
Using Mumford's Geometric Invariant Theory (GIT) from the work of Mumford~Fogarty~and~Kirwan~\cite{mum}, 
first Narasimhan~\&~Seshadri~\cite{nase}, and then Atiyah~\&~Bott~\cite{atbo}, build and charaterize this 
family of bundles by Morse theory and by their stability. Both classifications are equivalent. We use the 
one from stability here:

 \begin{Def} 
  Define the {\em moduli space of (polystable) vector bundles} $\sE \to \Sg$ as the quotient
   \[
   \sN(r,d) = \sA_{ps}^{0,1}(r,d) / \gG^{\bC}. 
   \]
   Respectively, define {\em the moduli space of stable vector bundles} as the quotient
   \[
   \sN_{s}(r,d) = \sA_{s}^{0,1}(r,d) / \gG^{\bC}
   \subseteq
   \sN(r,d). 
   \]
 \end{Def}
 
 \begin{Rmk} 
 The quotient
   $
   \sA^{0,1}(r,d) / \gG^{\bC}
   $
  is not Hausdorff so, we may lose a lot of interesting properties.
 \end{Rmk}
 
 Narasimhan~\&~Seshadri~\cite{nase} charaterize this variety, and then Atiyah~\&~Bott~\cite{atbo} conclude in their work that:
 
 \begin{Th}[9.~\cite{atbo}] 
  If $\GCD(r,d) = 1$, then $\sA_{ps}^{0,1} = \sA_{s}^{0,1}$
  and $\sN(r,d) = \sN_{s}(r,d)$ becomes a (compact) differential
  projective algebraic variety with dimension
  $$
  \dim_{\bC}\big(\sN(r,d)\big) = r^2 (g - 1) + 1.
  \QEDA
  $$
 \end{Th}

\subsection{Hitchin Construction} 
\label{ssec:3.2.}

Hitchin~\cite{hit2} works with the \emph{Yang-Mills self-duality equations} (SDE)
 \begin{equation} 
 \left\{
 \begin{array}{c c c}
  F_A + [\phi, \phi^{*}] & = & 0 \\
                         &   &   \\
  \bar{\partial}_A \phi  & = & 0
 \end{array}
 \right.
 \label{eq:YM} 
 \end{equation}
where $\phi \in \Omega^{1,0}\big(\Sg, \End(\sE)\big)$ is a complex auxiliary 
field and $F_A$ is the curvature of a connection $d_A$ which is compatible with 
the holomorphic structure of the bundle $E = (\sE, \bar{\partial}_A)$.  Hitchin 
calls such a $\phi$ as a \emph{Higgs~field}, because it shares a lot of the 
physical and gauge properties to those of the Higgs boson.

The set of solutions
 $$
 \beta(\sE) := 
 \{
 (\bar{\partial}_A, \phi)| \mathword{solution of}(\ref{eq:YM})
 \}
 \subseteq
 \sA^{0,1}(\sE) \times \Omega^{1,0}(\Sg,\End(\sE))
 $$
and the collection
 $$
 \beta_{ps}(r,d) := 
 \{
 \beta(\sE)| \sE \mathword{polystable,} \rk(\sE) = r, \deg(\sE) = d 
 \}
 $$
allow Hitchin to construct the Moduli space of solutions to SDE~(\ref{eq:YM})
  $$
  \sM^{YM}(r,d) = \beta_{ps}(r,d) / \gG^{\bC},
  $$
and
  $$
   \sM^{YM}_{s}(r,d) = \beta_{s}(r,d) / \gG^{\bC}
   \subseteq
   \sM^{YM}(r,d)
  $$
the moduli space of stable solutions to SDE~(\ref{eq:YM}).

\begin{Rmk}
 Recall that $\GCD(r,d) = 1$ implies $\sA_{ps}^{0,1} = \sA_{s}^{0,1}$
  and so $\sM^{YM}(r,d) = \sM^{YM}_{s}(r,d)$.
\end{Rmk}

Nitsure~\cite{nit} computes the dimension of this space in his work:

\begin{Th}[{Nitsure~(1991)}]
 The space $\sM(r,d)$ is a quasi--projective variety of complex dimension
 $$
 \dim_{\bC}(\sM(r,d)) = (r^2-1)(2g - 2).
 $$ 
\end{Th}

\subsection{Simpson Construction} 
\label{ssec:3.3.}

An alternative concept arises in the work of Simpson~\cite{sim2}:

\begin{Def}
 A \emph{Higgs bundle} over $\Sg$ is a pair $(E, \varphi)$ where 
 $E \to \Sg$ is a holomorphic vector bundle together with $\varphi$,
 an endomorphism of $E$ twisted by~$K = T^{*}(\Sg)\to \Sg$, the canonical 
 line bundle of the surface $\Sg$: $\varphi\: E \to E \ox K$. The field 
 $\varphi$ is what Simpson calls {\it Higgs field}.
\end{Def}

There is a condition of stability analogous to the one for vector bundles, 
but with reference just to subbundles preserved by the endomorphism $\varphi$:

\begin{Def}
 A subbundle $F \subset E$ is said to be \emph{$\varphi$-invariant} 
 if $\varphi(F) \subset F \ox K$. A Higgs bundle is said to be 
 \emph{semistable} [respectively, \emph{stable}] if 
 $\mu(F) \leq \mu(E)$ [resp., $\mu(F) < \mu(E)$] for 
 any nonzero $\varphi$-invariant subbundle $F \subseteq E$
 [resp., $F \subsetneq E$]. Finally, $(E,\varphi)$ is called
 \emph{polystable} if it is the direct sum of stable
 $\varphi$-invariant subbundles, all of the same slope.
\end{Def}

With this notion of stability in mind, Simpson~\cite{sim2} constructs
the Moduli space of Higgs bundles as the quotient
   $$
   \sM^{H}(r,d) = \{(E,\varphi)|\ E \mathword{polystable} \} / \gG^{\bC}
   $$
and the subspace
   $$
   \sM^{H}_{s}(r,d) = 
   \{(E,\varphi)|\ E \mathword{stable} \} / \gG^{\bC}
   \subseteq
   \sM^{H}(r,d)
   $$
of stable Higgs bundles.

\begin{Rmk}
  Again, $\GCD(r,d) = 1$ implies $\sM^{H}(r,d) = \sM^{H}_{s}(r,d)$. See 
  Simpson~\cite{sim2} for details.
\end{Rmk}

Finally, Simpson~\cite{sim2} concludes:

 \begin{Prop}[Prop.~1.5.~\cite{sim2}]
  There is a homeomorphism of topological spaces
  $$
   \sM^{H}(r,d)
   \cong
   \sM^{YM}(r,d).
   \QEDA
  $$
 \end{Prop}

\section{Recent Results} 
\label{sec:4.}

If we consider a fixed point $p \in X$ as a divisor $p \in \mathrm{Sym}^{1}(X)=X$, and $L_p$ the line bundle that corresponds to that divisor $p$, we get a complex of the form
\[
  E\xrightarrow{\quad \Phi^k \quad} E\otimes K \otimes L_p^{\otimes k}
\]
where $\Phi^k \in H^0(X,\mathrm{End}(E)\otimes K\otimes L_p^{\otimes k})$ is a \textit{Higgs field with poles of order $k$}. So, we call 
such a complex as a $k$-{\it Higgs bundle} and $\Phi^k$ as its $k$-{\it Higgs field}. As well as for a Higgs bundle, a $k$-Higgs bundle $(E,\Phi^k)$\ is stable (respectively 
semistable) if the slope of any $\Phi^k$-invariant subbundle of $E$ is strictly less (respectively less or equal) than the slope of $E:\ \mu(E)$.\ Finally, 
$(E,\Phi^k)$\ is called \textit{polystable} if $E$ is the direct sum of stable $\Phi^k$-invariant subbundles, all of the same slope.

Motivated by the results of Bradlow,~Garc\'ia--Prada,~Gothen~\cite{bgg2} and the work of Hausel~\cite{hau}, in terms of the homotopy groups of the moduli space of Higgs bundles and 
$k$-Higgs bundles respectively, the author has tried to improve the result of Hausel to higher rank. So far, here is one of the main results. Recall that here, only announce the result, no proof is given. The reader can find the proof in \cite{z-r2}.

\begin{Th}[Cor.~4.14.~\cite{z-r2}]
Suppose the rank is either $r = 2$ or $r = 3$, and $\GCD(r,d) = 1$. 
Then, for all $n$ exists $k_0$, depending on $n$, such that 
 $$
 \pi_j\left(\sM_{\Lambda}^{k}(r,d)\right) \xrightarrow{\quad \cong \quad} \pi_j\left(\sM_{\Lambda}^\infty(r,d)\right)
 $$
 for all $k \geq k_0$ and for all $j \leq n-1$.\QEDA
\end{Th}  

There is a lot of interesting recent results related with the homotopy of the moduli space of Higgs bundles, in terms of its stratifications (see~\cite{gzr}), and in terms of Morse Theory and Variations of Hodge Structures (see~\cite{z-r2}). Nevertheless, the introduction of a wide bunch of notation would be needed. The author encourages the reader to explore those results in the references above mentioned.

\renewcommand{\refname}{References}

\begin{flushright}
  \emph{Ronald A. Z\'u\~niga-Rojas}\\
  \small Centro de Investigaciones Matem\'aticas\\ 
  \small y Metamatem\'aticas CIMM\\
  \small Project: {\tt 820-B8-224}\\
  \small Escuela de Matem\'atica\\ 
  \small Universidad de Costa Rica UCR\\
  \small San Jos\'e 11501, Costa Rica\\
  \small e-mail: \texttt{ronald.zunigarojas@ucr.ac.cr}
\end{flushright}


\begin{thebibliography}{T-W}
\addcontentsline{toc}{section}{References}

\bibitem{atbo}
M. F. Atiyah and R. Bott,
Yang-Mills equations over Riemann surfaces,
\emph{Phil. Trans. Roy. Soc. London}, A \textbf{308} (1982), 523--615.


\bibitem{bgg}
S.B. Bradlow, O. Garc\'ia-Prada and P.B. Gothen,
What is a Higgs Bundle?,
\emph{Notices of the AMS}, \textbf{54}, No.~8 (2007), 980--981.


\bibitem{bgg2} 
S. B. Bradlow, O. Garc\'ia-Prada and P. B. Gothen,
Homotopy groups of moduli spaces of representations,
\emph{Topology} \textbf{47} (2008), 203--224.


\bibitem{gzr}
P. B. Gothen and R.~A.~Z\'u\~niga-Rojas,
Stratifications on the moduli space of Higgs bundles,
\emph{Portugaliae Mathematica, EMS}, \textbf{74}, (2017), 127--148.


\bibitem{hau} 
T. Hausel,
\emph{Geometry of Higgs bundles}, 
Ph.\,D. thesis, University of Cambridge, UK, 1998.


\bibitem{hit1} 
N.J.~Hitchin, 
Gauge Theory on Riemann Surfaces, 
\emph{Proceedings of the First College of Riemann Surfaces} held in Trieste, Italy $(1987)$ 99-118.


\bibitem{hit2} 
N. J. Hitchin, 
The self-duality equations on a Riemann surface,
\emph{Proc. London Math. Soc.}, \textbf{55} (1987), 59--126.


\bibitem{kob}
S.~Kobayashi, 
\emph{Differential Geometry of Complex Vector Bundles},
Publications of the Mathematical Society of Japan, Vol.~{\bf 15}, Iwanami~Shoten~Publ. and Princeton~Univ.~Press, 1987.


\bibitem{lute}
M.~L\"ubke~and~A.~Teleman,
\emph{The Kobayashi-Hitchin Correspondence},
World Scientific Publishing Co., 1995.


\bibitem{mum}
D. Mumford, J. Fogarty and F. Kirwan,
\emph{Geometric Invariant Theory},
Springer, 1994.


\bibitem{nase} 
M.S. Narasimhan and C.S. Seshadri, 
Stable and Unitary Vector Bundles on a Compact Riemann Surface,
\emph{The Annals of Mathematics}, Second Series Vol. \textbf{82} (1965), pp. 540-567.


\bibitem{nit} 
Nitsure, N., 
Moduli Space of Semistable Pairs on a Curve,
\emph{Proc. London Math. Soc.}, (3), \textbf{62}, (1991), pp. 275-300.


\bibitem{sim2} 
C. T. Simpson, 
Higgs bundles and local systems,
\emph{Publ. Math. IH\'ES} \textbf{75} (1992), 5--95.


\bibitem{yami}
C.~N.~Yang and R.~L.~Mills,
Conservation of isotopic spin and isotopic gauge invariance,
{\em Phys.~Rev.} {\bf 96}, $(1954)$, 191--195.


\bibitem{z-r0}
R.~A.~Z\'u\~niga-Rojas,
\emph{Homotopy groups of the moduli space of Higgs bundles},
Ph.\,D. thesis, Universidade do Porto, Portugal, 2015.


\bibitem{z-r1}
R.~A.~Z\'u\~niga-Rojas,
On the Cohomology of The Moduli Space of sigma-Stable Triples and (1,2)-VHS, to appear.
\texttt{arXiv:1803.01936 [math.AG]}.


\bibitem{z-r2}
R.~A.~Z\'u\~niga-Rojas,
Stabilization of the Homotopy Groups of The Moduli Space of k-Higgs Bundles, 
\emph{Revista Colombiana de Matem\'aticas},  {\bf 52} (2018) 1, 9-31. 

\end{thebibliography}
\end{document}